\documentclass[12pt]{article}%
\usepackage{amsmath}
\usepackage{amsfonts}
\usepackage{amssymb}
\usepackage{graphicx}%
\providecommand{\U}[1]{\protect\rule{.1in}{.1in}}
\newtheorem{theorem}{Theorem}

\newtheorem{corollary}{Corollary}

\newtheorem{definition}{Definition}

\newtheorem{proposition}{Proposition}
\newtheorem{remark}{Remark}

\begin{document}

\title{On accompanying measures and asymptotic expansions in limit theorems for
maximum of random variables.}
\author{V. I. Piterbarg\thanks{Lomonosov Moscow state university, Laboratory of
Stochastic Analysis and its Applications of Higher School of Economics,
Moscow; Scientific Research Institute of System Analysis of Russia;
\lowercase{e-mail:
piter@mach.math.msu.su}}, Yu. A. Scherbakova\thanks{Lomonosov Moscow state
university, 119991, Moscow}}
\maketitle

\begin{abstract}
\footnotetext[1]{Partially funded by RFBR and CNRS, project number
20-51-15001} A sequence of accompanying laws is suggested in the limit theorem
of B. V. Gnedenko for maximums of independent random variables belonging to
maximum domain of attraction of the Gumbel distribution. It is shown that this
sequence gives an exponential power rate of convergence whereas the Gumbel
distribution gives only a logarithmic rate. As examples, classes of Weibull
and log-Weibull type distributions are considered in details. A scale for the
Gumbel maximum domain of attraction is suggested as a continuation of the
considered two classes.

\emph{{Keywords: Gnedenko-Fisher-Tippet theorem, rate of convergence,
correction terms, accompanying law } }

\end{abstract}

\section{  Introduction}

Let $X_{1},...,X_{n},...,$ be independent identically distributed random
variables having distribution function $F(x).$ Denote $M_{n}:=\max
(X_{i},i=1,...,n).$ By Gnedenko-Fisher-Tippet theorem, see \cite{Gnedenko1},
and also \cite{deHaanFer}, \cite{EKM}, if for some positive $a_{n}$ and real
$b_{n}$ there exists a non-degenerate limit of $P(M_{n}\leq a_{n}x+b_{n})$ as
$n\rightarrow\infty,$ then the limit distribution function belongs to one of
three types, that is Fr\`{e}chet, Weibull, and Gumbel types. We deal in this
note with the third one, that is, with distribution functions $F$ for which
for any $x,$
\begin{equation}
\lim_{n\rightarrow\infty}P(M_{n}\leq a_{n}x+b_{n})=\Lambda(x):=\exp(-e^{-x}).
\label{Gumbel_limit}%
\end{equation}
One says in this case that the distribution belongs to Gumbel max domain of
attraction, $F\in MDA(\Lambda).$ Moreover, we restrict ourselves by the case
of infinite right end point of $F,$ that is, $F(x)<1$ for any positive $x.$ It
will be seen that our approach can be applied to any distributions from the
three above mentioned domains of maximum attraction. We take this domain
because of it is extremely wide comparing with other two. Remark that many
applications require a division of this domain into reasonable in some sense
sub-domains. For example, this domain contains distributions with tails
$1-F(x),$ roughly (logarithmic) equivalent to Weibull tails, that is,
$\log(1-F(x))\sim-Cx^{p},$ $C,p>0,$ as $x\rightarrow\infty$, or log-Weibull
tails, that is, $\log(1-F(x))\sim-C\log^{p}x,$ $C>0,p>1,$ wherein one can take
regularly varying on infinity functions instead of degrees in the right-hands
in the above relations, and many other distributions with heavier and lighter
tails. for comparison, recall that the Frech\'{e}t maximum domain of
attraction contains only distributions with regularly varying tails on
infinity, and Weibull maximum domain of attraction contains only distributions
with regularly varying tails at the right (finite) end of the distribution. We
discuss here an idea of \textit{scaling} the distributions from $MDA(\Lambda)$
with mentioned above first two "scale points", that is, generalized
Weibull-like and generalized log-Weibull-like distributions, see definitions below.

By \cite{BdH}, the distributions from $MDA(\Lambda)$ can be described using
von Mises function, that is, $F\in MDA(\Lambda)$ if and only if for some
$x_{0}\geq0$,
\begin{equation}
1-F(x)=c(x)\exp\left\{  -\int_{x_{0}}^{x}\frac{1}{f(t)}dt\right\}  ,\ \ x\geq
x_{0}, \label{fon Mises0}%
\end{equation}
with $f(x),$ positive and absolutely continuous on $[x_{0},$ $\infty)$ with
density $f^{\prime}(x),$ $f^{\prime}(x)\rightarrow0,$ and $c(x)\rightarrow
c>0,$ as $x\rightarrow\infty$.

Often more flexible form of this assertion is convenient: $F\in MDA(\Lambda)$
if and only if for some $x_{0}\geq0$,
\begin{equation}
1-F(x)=c(x)\exp\left\{  -\int_{x_{0}}^{x}\frac{g(t)}{f(t)}dt\right\}  ,
\label{fon MisesG}%
\end{equation}
with the same properties of $f(x),$ and $c(x)$ and $g(x)\rightarrow1,$ as
$x\rightarrow\infty.$ For both the representations one can take%
\begin{equation}
b_{n}=F^{\leftarrow}(1-n^{-1}),\ \ \ a_{n}=f(b_{n}). \label{constants}%
\end{equation}
See \cite{EKM}, \cite{Resnick}. Moreover, if (\ref{Gumbel_limit}) is fulfilled
with $b_{n}$ as in (\ref{constants}) and some sequence of positive $\tilde
{a}_{n}$ then for $a_{n}$ as in (\ref{constants}), $a_{n}/\tilde{a}%
_{n}\rightarrow1$ as $n\rightarrow\infty$ and there exists a representation
(\ref{fon MisesG}) with other $\tilde{g}$ and $\tilde{f}$ such that $\tilde
{a}_{n}=\tilde{f}(b_{n}).$

There is a wide bibliography on the quality of convergence in
(\ref{Gumbel_limit}). We note here two main directions of studies. First one
is related with restrictions on the tail behavior of $F$ at infinity. First of
all it is the second order condition and higher orders regular behaviors of
the tails of $F,$ see \cite{deHaanResnick}, \cite{LJPeng}, \cite{Alves et
all}, \cite{deHaanFer}, \cite{WangCheng}, and references therein. Definitions
and comments see also in section 3 below.

The second direction related to concrete expressions of the distributions or
families of distributions, such as Gaussian, Gaussian like, Weibull,
Weibull-like, log-Weibull-like, so on. Here is also a wide bibliography, beginning with
\cite{Hall}, \cite{HallWellner}, \cite{PengNadara}, see also monographs
\cite{deHaanFer} and \cite{Resnick}.

Our study belongs rather to the second direction, we use von Mises structure
(\ref{fon Mises0}) and (\ref{fon MisesG}), but we suggest another approach: do
not investigate quality of Gumbel approximation but first look for better
approximations. Notice that this is very common approach in the study of
quality approximation in Central Limit Theorem. This is Chebyshev-Hermite
polynomials approximation, other types of approximation, other types of
accompanying laws and charges (signed measures), see \cite{Feller}, \cite{GK},
\cite{IbragimovPresman}, \cite{petrov}, \cite{senatov-book}, \cite{senatov}.
In connection with, in 2002, one of the authors discussed with Laurens de Haan
the following result on Gaussian smooth stationary processes, only recently it
is published in \cite{lect}, we then have agreed that such approach is
interesting and promising.

\begin{theorem}
\label{main} Let $X(t),$ $t\in\mathbb{R},$ be a twice differentiable in square
mean Gaussian stationary process with $EX(t)=0,$ $EX^{2}(t)=1,$ $EX^{\prime
}(t)^{2}=1.$ Assume that
\[
\int_{0}^{\infty}|r(t)|^{a}dt<\infty
\]
holds for its covariance function $r$ and some $a>0$. Denote $l_{T}%
=\sqrt{2\log\frac{T}{2\pi}}.$ Denote also
\begin{equation}
A_{T}(x)=\left\{
\begin{array}
[c]{c}%
e^{-e^{-x-x^{2}/2l_{T}^{2}}},\ x\geq-l_{T}^{3/2},\\
0,\ x<-l_{T}^{3/2},
\end{array}
\right.  \label{A_T}%
\end{equation}
$T>0$. Then

\begin{enumerate}
\item For some $\gamma>0$,
\[
P\left(  \max_{t\in\lbrack0,T]}X(t)\leq l_{T}+\frac{x}{l_{T}}\right)
-A_{T}(x)=O(T^{-\gamma}),\ \ \ T\rightarrow\infty
\]
uniformly in $x\in\mathbb{R}.$

\item Moreover,
\[
l_{T}^{2}\left(  P\left(  \max_{t\in\lbrack0,T]}X(t)\leq l_{T}+\frac{x}{l_{T}%
}\right)  -e^{-e^{-x}}\right)  \rightarrow\frac{1}{2}e^{-e^{-x}}e^{-x}x^{2},
\]
as $T\rightarrow\infty$, uniformly in $x\in\mathbb{R}.$
\end{enumerate}
\end{theorem}

It follows from the second statement that the rate of convergence of the
distribution of the maximum to the Gumbel distribution is logarithmic. It also
gives the second term of the asymptotic expansion for the probability. The
first statement gives the sequence of approximating functions that approaches
the maximum distribution with the power rate. There are several results
related to accompanying laws, see \cite{deHaanResnick}, \cite{ChengJiang},
\cite{DdHL}, \cite{DdHL1}, where the mentioned above second order condition is
exploited. Many of the results are presented in celebrated monographs
\cite{Resnick} and \cite{deHaanFer}.

In the following section, an asymptotic expansion is derived in the theorem on
convergence to Gumbel distribution and speed of convergence is studied. It
turns out that the quality of convergence is as a rule logarithmic. Further,
by analogy with corresponding results on Central Limit Theorem, a sequence of
accompanying charges (signed measures) is introduced. This sequence gives a
power rate of convergence. Comparing our results with the known ones is
considered in Section 3. Section 4 contains examples for distributions with
smooth tails. Some of the examples of particular distributions were subjects
of student works in faculty of mechanics and mathematics of Lomonosov Moscow
state university, including one of the authors of the present work. We thank
Viktoria Maier, Ignat Melnilov, Viktor Troshin, Kirill Lisakov for theirs
help. Finally, we discuss in Section 5 a scale for distributions from
$MDA(\Lambda)$ and some related problems.

\section{Asymptotic expansions and accompanying measures}

In contrast to similar problems related to the central limit theorem,
\cite{petrov}, \cite{IbragimovPresman}, \cite{senatov}, \cite{senatov-book},
the construction of asymptotic expansions and accompanying measures in limit
theorems for maximums is much simpler, in a certain sense even trivial.
Similar evaluations based on Taylor expansions one can find in many works on
extreme distributions, see, for example, \cite{novak} and references given in
Section 3. Nevertheless, we present these evaluations here since they are some
background for similar calculations for particular distributions and
distribution classes.

Assume (\ref{fon Mises0}). By (\ref{constants}),
\begin{equation}
\int_{x_{0}}^{b_{n}}\frac{1}{f(t)}dt=\log(nc(b_{n})). \label{d_n}%
\end{equation}
Further,
\begin{equation}
G_{n}(x):=F^{n}(a_{n}x+b_{n})=\left(  1-c(a_{n}x+b_{n})e^{-\int_{x_{0}}%
^{a_{n}x+b_{n}}\frac{1}{f(t)}dt}\right)  ^{n}. \label{G_n def}%
\end{equation}
Taking logarithm, and denoting
\[
g_{n}(x):=\int_{x_{0}}^{a_{n}x+b_{n}}\frac{1}{f(t)}dt-\log c(a_{n}x+b_{n}),\
\]
we have after easy calculations using Taylor,
\begin{equation}
\log G_{n}(x)=-ne^{-g_{n}(x)}\sum_{k=0}^{\infty}\frac{1}{k+1}e^{-kg_{n}(x)}.
\label{logG}%
\end{equation}
Using (\ref{d_n}),%
\begin{align*}
g_{n}(x)  &  =\int_{x_{0}}^{b_{n}}\frac{1}{f(t)}dt\\
&  +\int_{b_{n}}^{a_{n}x+b_{n}}\frac{1}{f(t)}dt-\log c(a_{n}x+b_{n})\\
&  =\log n+\int_{b_{n}}^{a_{n}x+b_{n}}\frac{1}{f(t)}dt-\log\frac
{c(a_{n}x+b_{n})}{c(b_{n})}=:\log n+\gamma_{n}(x).
\end{align*}
Now (\ref{logG}) can be written as
\[
\log G_{n}(x)=-e^{-\gamma_{n}(x)}\sum_{k=0}^{\infty}\frac{1}{(k+1)n^{k}%
}e^{-k\gamma_{n}(x)}.
\]
Taking out the first summand in the sum and passing to exponents, we get,
that
\begin{equation}
G_{n}(x)=\exp\left(  -e^{-\gamma_{n}(x)}\right)  \exp\left(  -\frac{1}%
{n}\Sigma(x)\right)  , \label{G_n(x)}%
\end{equation}
with
\begin{equation}
\Sigma(x)=\sum_{k=0}^{\infty}\frac{\exp(-(k+2)\gamma_{n}(x))}{(k+2)n^{k}}.
\label{Sigma(x)}%
\end{equation}
From here and (\ref{G_n def}), using Taylor, we get that ,
\begin{align}
P(M_{n}  &  \leq a_{n}x+b_{n})=\exp\left(  -e^{-\gamma_{n}(x)}\right)
\nonumber\\
&  +\frac{1}{n}\exp\left(  -e{}^{-\gamma_{n}(x)}\right)  \sum_{k=0}^{\infty
}\frac{(-1)^{k+1}\Sigma^{k+1}(x)}{(k+1)!n^{k}}. \label{expansion1_2}%
\end{align}
Writing summands in (\ref{Sigma(x)}) as $\exp(-(k+2)\gamma_{n}(x)-\log
(k+2)-k\log n),$ one can see that
\[
\sup_{x:\gamma_{n}(x)\geq-\log n}|\Sigma(x)|<\infty,
\]
hence the same is valid for the sum the second term on the right hand part of
(\ref{expansion1_2}). Introduce therefore the \emph{accompanying sequence},
\begin{equation}
B_{n}(x)=\left\{
\begin{array}
[c]{c}%
e^{-e^{-\gamma_{_{n}}(x)}},\ \gamma_{n}(x)\geq-\log n,\\
0,\ \gamma_{n}(x)<-\log n.
\end{array}
\right.  \label{B_n(x)}%
\end{equation}

Turn now to $\gamma_{n}(x)$. Notice first that for any $x,$ by
(\ref{Gumbel_limit}) and (\ref{expansion1_2}), $\gamma_{n}(x)\rightarrow x$ as
$n\rightarrow\infty.$ Further, since $a_{n}=f(b_{n}),$ rewrite the expression
for $\gamma_{n}(x)$ as%
\begin{equation}
\gamma_{n}(x)=\int_{b_{n}}^{a_{n}x+b_{n}}\left(  \frac{1}{f(t)}-\frac
{1}{f(b_{n})}\right)  dt-\log\frac{c(a_{n}x+b_{n})}{c(b_{n})}+x.
\label{gamma_n}%
\end{equation}
One can easily see that, assuming (\ref{fon MisesG}) instead of
(\ref{fon Mises0}), the expression for $a_{n}$ is changed on $a_{n}%
=f(b_{n})/g(b_{n})$, and the expression for $\gamma_{n}(x)$ is changed on
\begin{equation}
\gamma_{n}(x)=\int_{b_{n}}^{a_{n}x+b_{n}}\left(  \frac{g(t)}{f(t)}%
-\frac{g(b_{n})}{f(b_{n})}\right)  dt-\log\frac{c(a_{n}x+b_{n})}{c(b_{n})}+x.
\label{gamma_ng}%
\end{equation}
Changing in (\ref{gamma_n}) variables $t=b_{n}+v/b_{n}$ and using
$a_{n}=f(b_{n}),$ we have after some algebra,%
\begin{equation}
\gamma_{n}(x)=\int_{0}^{x}\left(  \frac{f(b_{n})}{f(b_{n}+a_{n}v)}-1\right)
dv-\log\frac{c(a_{n}x+b_{n})}{c(b_{n})}+x; \label{gamma1}%
\end{equation}
and for general $g,$
\begin{equation}
\gamma_{n}(x)=\int_{0}^{x}\left(  \frac{g(b_{n}+a_{n}v)f(b_{n})}{f(b_{n}%
+a_{n}v)g(b_{n})}-1\right)  dv-\log\frac{c(a_{n}x+b_{n})}{c(b_{n})}+x.
\label{gamma1g}%
\end{equation}

Finally we get the following.

\begin{theorem}
\label{expansion} Let (\ref{fon Mises0}) or (\ref{fon MisesG}) be fulfilled
for a distribution function $F(x).$ Let $X_{1},X_{2},...,$ be i.i.d. random
variables with distribution function $F,$ and $M_{n}:=\max(X_{1},,...X_{n}).$
Then for $a_{n},b_{n}$ defined by (\ref{constants}),
\[
P(M_{n}\leq a_{n}x+b_{n})-B_{n}(x)=O\left(  1/n\right)  ,\ \text{as\ }%
n\rightarrow\infty,
\]
 uniformly in $x\in \mathbb{R}$, where $B_{n}(x)$ is defined by (\ref{B_n(x)}) with
$\gamma_{n}(x)$ defined by (\ref{gamma_n}) or (\ref{gamma_ng}) correspondingly
to the representations of $F$. Moreover, for all $x$, $B_{n}(x)\rightarrow
\Lambda(x)$ as $n\rightarrow\infty$.
\end{theorem}

That is, the sequence $B_{n}(x)$ is \emph{ a natural sequence of accompanying
charges (laws)} in Gnedenko Limit Theorem. which gives the power rate of
convergence to Gumbel distribution. Now one can see that the question on the
rate of convergence to Gumbel ditribution depends on how fast $\gamma_{n}(x)$
tends to $x$ as $n\rightarrow\infty.$ Definitely, it depends on detailed tail
behavior of $F(x).$ Generally, we may formulate an analogous of the second
statement of Theorem \ref{main} as follows.

\begin{corollary}
\label{rate of convergence} In Theorem \ref{expansion} conditions, for any
$x$,
\[
P(M_{n}\leq a_{n}x+b_{n})-\exp\left(  -e^{-x}\right)  =\exp\left(
-e^{-x}\right)  e^{-x}(\gamma_{n}(x)-x)(1+o(1))+O(n^{-1})
\]
as $n\rightarrow\infty$.
\end{corollary}

Indeed, it follows from Taylor expansion of $G_{n}(x).$ We point out that, in
our opinion, a formulation of complete analogous of Theorem \ref{main}
statement 2 is possible only under sufficiently concrete tail behavior
description. In Theorem \ref{main} the behavior is given exactly. Remark that
for exponential distribution, $f(x),$ $c(x)$ and $g(x)$ are constants, hence
$\gamma_{n}(x)\equiv x,$ that is the rate of convergence in
(\ref{Gumbel_limit}) is $O(n^{-1}),$ $n\rightarrow\infty,$ the exact
expression can be derived, using for example, (\ref{expansion1_2}). In Section
4 we prove several other refinements for distributions with different types 
of theirs tail behaviors, see Propositions  \ref{rate_weibull-like}, \ref{weibull_like_standad_rate} and \ref{log-Weib_prop}.

Thus in order to study the rate of convergence in (\ref{Gumbel_limit}) one has
first to study the behavior of $\gamma_{n}(x)-x$ as $n\rightarrow\infty$ and
in dependence of $x.$ Then one should compare the asymptotic behavior of the
difference with the behavior of the residual which is equal to $O(1/n)$.
Notice that the two first summand in (\ref{gamma_ng}) both can play the main
role in the rate of convergence.

Now give another obvious expression for $\gamma_{n}(x).$

\begin{proposition}
\label{expr_for_gamma} Let (\ref{Gumbel_limit}) be fulfilled, then for the
corresponding representation (\ref{fon MisesG}), with $a_{n}=f(b_{n}),$
\begin{equation}
\gamma_{n}(x)=-\log\frac{1-F(b_{n}+a_{n}x)}{1-F(b_{n})}=\log\frac
{1}{n(1-F(b_{n}+a_{n}x)}. \label{expression2_gamma}%
\end{equation}

\end{proposition}

Indeed, it follows from the above that
\[
\gamma_{n}(x)=\int_{x_{0}}^{a_{n}x+b_{n}}\frac{g(t)dt}{f(t)}-\int_{x_{0}%
}^{b_{n}}\frac{g(t)dt}{f(t)}-\log c(a_{n}x+b_{n})+\log c(b_{n}),
\]
then use (\ref{fon MisesG}).

\section{Relation to some known results.}

As it was already mentioned, one of main approaches to studies the rate of
convergence in limit theorems for maximums is introducing additional
conditions on behavior of tails pf distribution functions. Here we consider
how this conditions relate to the behavior of $\gamma_{n}(x).$ Notice that we
consider $F\in MDA(\Lambda)$ with $F(x)<1$ for all $x.$

\begin{definition}
\textbf{(Second order condition for functions from } $MDA(\Lambda)$). There
exists a function $A(n)$ of constant sign which tends to zero as
$n\rightarrow\infty$ and such that there exists the limit
\begin{equation}
\lim_{n\rightarrow\infty}\frac{e^{-\gamma_{n}(x)}-e^{-x}}{A(n)}=H(x),
\label{second order condition}%
\end{equation}
and $H(x)$ is not identicaly equals neither zero nor infinity.
\end{definition}

This condition is introduced by L. de Haan, \cite{deHaan1}, in some other
terms. We give equivalent formulation, based on Theorem 2.3.8,
\cite{deHaanFer}, and on represenrarion (\ref{expression2_gamma}) of
$\gamma_{n}(x).$ From this definition immediately follows, see for example,
\cite{deHaanFer}, that $A(n)$ regularly varies on infinity with non positive
index $\rho\leq0.$ It is also known that for the considered here case of
convergence to Gumbel distribution,
\begin{equation}
H(x)=\frac{1}{\rho}\left(  \frac{x^{\rho}-1}{\rho}-\log x\right)
,\ \text{åñëè\ }\rho<0,
\label{second order conditionA}%
\end{equation}
and
\[
H(x)=\frac{1}{2}\log^{2}x\ \ \text{åñëè\ }\rho=0.
\]
Using above mentioned Theorem 2.3.8 and Proposition \ref{expr_for_gamma}, we get
the following.

\begin{corollary}
\label{second order} Let conditions of Theorem \ref{expansion} and
(\ref{second order condition}) be fulfilled. Then
\begin{align*}
P(M_{n}  &  \leq a_{n}x+b_{n})=\exp\left\{  -e^{-x}-A(n)H(x)(1+o(1)\right\} \\
&  \times\exp\left(  -\frac{1}{n}{}\sum_{k=0}^{\infty}\frac{1}{(k+2)n^{k}%
}\left(  \frac{1-F(a_{n}x+b_{n})}{1-F(b_{n})}\right)  ^{k+2}\right)  .
\end{align*}

\end{corollary}

Remark that if $\rho<-1$ the second exponent gives the main contribution in
the rate of convergence; if $\rho=-1$, one should knows the behavior of
$A(n)=n^{-1}\ell(n)$ more exactly, that is, the behavior of the corresponding
slowly varied function $\ell(n)$; finally, if $\rho>-1,$ then the second
summand in the first exponent gives the main contribution.

Similar calculations can be performed also for the $n$th order condition, see
\cite{WangCheng}.

Let us give another result on the rate of convergence and accompanying laws,
it is interesting to compared it with Theorem \ref{main}.

\begin{theorem}
\label{Theorem 2.1}(\textbf{Theorem 2.1, \cite{DdHL}) } Let
(\ref{second order condition}) be fulfilled ans let $\rho<0,$ see
(\ref{second order conditionA}). Take $b_{n}=F^{\leftarrow}\left(
e^{-1/n}\right)  $ and correspondingly $a_{n}=f(b_{n})$. Then for any
$\varepsilon>0$ the relation
\[
\sup_{x}e^{(1-\varepsilon)x}\left\vert \frac{F^{n}(a_{n}x+b_{n})-\exp
(-e^{-x})}{A(n)}+\frac{1}{\rho}e^{-x+\rho x}e^{-e^{-x}}\right\vert
\rightarrow0
\]
takes place as $n\rightarrow\infty.$
\end{theorem}

\begin{remark}
\label{question2} This theoren, as well as Theorem 5.3.3, \cite{deHaanFer},
can be obtained immediately from Theorem \ref{expansion}. It is follow from
proofs of the mentioned theorems.
\end{remark}

\begin{remark}
It is interesting to use Theorem \ref{expansion} in studying large deviations
probabilities in Gnedenko Limit Theorem. For example, Corollary 2.1,
\cite{DdHL} and Theorem 5.3.12, \cite{deHaanFer}, under appropriate
restrictions, can be obtained from suggested here asymptotical expansions.
Notice also, that in \cite{Resnick}, similar asymptotic expansions are used
for this purpose.
\end{remark}

\section{Distributions with absolutely continuous tails \label{smooth}}

Assuming that $F$ is eventually, for large $x,$ absolutely continuous, hence
$c(x),$ (\ref{fon Mises0}, \ref{fon MisesG}) is eventually absolutely
continuous, and since it is additionally eventually strictly positive, there
exists $x_{0}$ such that for all $x\geq x_{0},$
\[
\log c(x)-\log c(x_{0})=\int_{x_{0}}^{x}\frac{c^{\prime}(t)dt}{c(t)}.
\]
Hence the representation (\ref{fon Mises0}) can be easily transformed to
\[
1-F(x)=c\exp\left\{  -\int_{x_{0}}^{x}\frac{c(x)-c^{\prime}(x)f(x)}%
{c(x)f(t)}dt\right\}  ,\ \ x\geq x_{0},
\]
where we have re-defined $c(x)$ on $cc(x),$ and change $x_{0}$ to have
positive and absolutely continious $c^{\prime}(x)$ for all $x\geq x_{0},$
with
\begin{equation}
f(x)c^{\prime}(x)\rightarrow0,\text{ and }f^{2}(x)c^{\prime\prime
}(x)\rightarrow0\ \text{as }x\rightarrow\infty. \label{Fsmooth}%
\end{equation}
In this conditions, the function
\[
\tilde{f}(x):=\frac{f(x)}{1-\frac{c^{\prime}(x)}{c(x)}f(x)}%
\]
satisfies the same conditions as $f,$ so that we can write%
\begin{equation}
1-F(x)=c\exp\left\{  -\int_{x_{0}}^{x}\frac{1}{\tilde{f}(t)}dt\right\}
,\ \ x\geq x_{0}. \label{fon Mises0smooth}%
\end{equation}
That is, \emph{a ditribution function from the Gumbel domain of maximal
attraction which has eventually sufficiently smooth tail, with introduced
above restrictions on }$c(t)$\emph{ is a von Mises function.} Having this
\emph{strengthened von Mises representation}, we also have a corresponding
shortened representation for $\gamma_{n}(x),$%
\begin{equation}
\gamma_{n}(x)=\int_{0}^{x}\left(  \frac{\tilde{f}(b_{n})}{\tilde{f}(b_{n}%
+v)}-1\right)  dv+x. \label{gamma_nsmooth}%
\end{equation}
We saw above that the behavior of $c(x)$ may give main contribution to the
rate of convergence. On the contrary, in the introduced in this Section
conditions, the influences of all the functions $c(x),$ $g(x)$ and $f(x)$ are
aggregated by the function $\tilde{f}(x).$ So, let us consider the behavior of
function
\begin{equation}
\gamma(t;x):=\int_{0}^{x}\left(  \frac{f(t)}{f(t+v)}-1\right)
dv\ \ \label{gammatx}%
\end{equation}
as $t\rightarrow\infty.$ By (\ref{fon Mises0}, \ref{Gumbel_limit}) and
Proposition \ref{expansion}, for any $x,$ $\gamma(t;x)\rightarrow0.$ Moreover,
since now $f$ is also eventually differentiable,
\[
\frac{f(t+v)}{f(t)}-1=\int_{0}^{v}\frac{f^{\prime}(t+s)}{f(t)}ds,
\]
hence, for almost all $s,$
\[
\frac{f^{\prime}(t+s)}{f(t)}\rightarrow0
\]
as $t\rightarrow\infty.$ Therefore, for smooth tails, having (\ref{Fsmooth}),
we may investigate the rate of convergence only in terms of $f.$ Below we
consider two important classes of distributions which are subsets of
$MDA(\Lambda),$ namely, ditributions with \emph{Weibull-like }and
\emph{log-Weibull-like }tails. Then we suggest a corresponding scale for
$MDA(\Lambda).$

Notice finally that since only ultimately behavior of the tail is in the frame
of our consideration, we may modify $f$ and $c$ on any bounded fixed interval,
in dependence of questions under consideration. For example we may put
$x_{0}=0$ or $x_{0}=1,$ changing correspondingly the functions in the von
Mises representations.

\subsection{Generalized Weibull-like distributions}

Below on we assume that $F$ is absolutely continuous. For flexible von Mises
representation (\ref{fon MisesG}), one can choose appropriately $g(t)$ on a
finite interval, say, $[0,a],$ to have $c(x)\equiv1,$ $x_{0}=1.$ Now take
$f(t)=Ct^{1-p}$ with $p,C>0$. That is, for all $x\geq1$,
\begin{equation}
1-F(x)=\exp\left(  -\int_{1}^{x}\frac{g(t)dt}{Ct^{1-p}}\right)  .
\label{weib_like}%
\end{equation}
Such the distributions we call generalized Weibull-like distributions.

Let us find norming sequences $a_{n}$ and $b_{n},$taking in mind that in
virtue of covvergence to types theorem, see \cite{GK}, \cite{petrov},
\cite{EKM}, the limit in (\ref{Gumbel_limit}), after changing constants
$(a_{n},b_{n})$ on $(\tilde{a}_{n},\tilde{b}_{n}),$ belongs to the same
(Gumbel) type, if and only if
\begin{equation}
\frac{a_{n}}{\tilde{a}_{n}}\rightarrow1,\ \frac{b_{n}-\tilde{b}_{n}}{a_{n}%
}\rightarrow0. \label{types}%
\end{equation}

Denote $\alpha(x):=g(x)-1.$ Integrating, we have,
\begin{equation}
1-F(x)=\exp\left(  \int_{1}^{x}\frac{\alpha(t)dt}{Ct^{1-p}}\right)
\exp\left(  -\frac{x^{p}}{Cp}\right)  . \label{weib_like1}%
\end{equation}
Since $\alpha(t)\rightarrow0$ as $t\rightarrow\infty,$
\[
\int_{1}^{x}\frac{\alpha(t)dt}{Ct^{1-p}}=o(x^{p})
\]
as $x\rightarrow\infty.$ An equation for $b_{n}$ is
\[
\int_{1}^{x}\frac{g(t)dt}{t^{1-p}}=C\log n.
\]
Passing to $\alpha(t),$ integrating in parts, and denoting $y=x^{p},$ we come
to the following equation,
\begin{equation}
y+\int_{1}^{y^{1/p}}\frac{p\alpha(t)dt}{t^{1-p}}=Cp\log n+p. \label{iter1}%
\end{equation}
From here we have, $y=Cp\log n(1+o(1))$ as $n\rightarrow\infty.$ Now we apply
the asymptotical iteration as $n\rightarrow\infty$, see for example
\cite{olver}. Denote $u:=Cp\log n,$ write for convenience $y=u(1+\varepsilon
(y,u))^{p}$, and find $\varepsilon(y,u).$ Equation (\ref{iter1}) becomes as
following,
\begin{align}
y  &  =u-\int_{1}^{u^{1/p}(1+\varepsilon(y,u)}\frac{p\alpha(t)dt}{t^{1-p}%
}=u-\int_{1}^{u^{1/p}}\frac{p\alpha(t)dt}{t^{1-p}}-\int_{u^{1/p}}%
^{u^{1/p}(1+\varepsilon(y,u))}\frac{p\alpha(t)dt}{t^{1-p}}\nonumber\\
&  =u-\int_{1}^{u^{1/p}}\frac{p\alpha(t)dt}{t^{1-p}}-u^{-1}\int_{1}%
^{1+\varepsilon(y,u)}\frac{p\alpha(u^{-1/p}s)ds}{s^{1-p}}. \label{iter2}%
\end{align}
Since we may assume that $\alpha(1)>0,$ the latter integral is of order
$\varepsilon(y,u).$ Now, in order to find $\varepsilon(y,u),$ one can again
input this expression for $y$ into (\ref{iter1}), so on, repeating this
several times.

Finally we get from equation (\ref{iter2}) that
\begin{equation}
b_{n}^{p}=Cp\log n+p-\int_{1}^{(Cp\log n)^{1/p}}\frac{p\alpha(t)dt}{t^{1-p}%
}-\frac{\varepsilon^{\prime}(n,b_{n})}{Cp\log n}, \label{b_n}%
\end{equation}
with
\[
\varepsilon^{\prime}(n,b_{n})=\int_{1}^{1+\varepsilon(b_{n}^{p},Cp\log
n)}\frac{p\alpha((Cp\log n)^{-1/p}s)ds}{s^{1-p}}.
\]
Notice that since $\alpha(t)$ vanishes as $t$ increases, the second term on
the right in (\ref{b_n}) is infinitely smaller the first one. For example,
depending on value of $p$ and on rate of tending $\alpha(t)$ to zero, one can
take
\[
b_{n}=(Cp\log n)^{1/p}-(Cp\log n)^{1/p-1}\int_{1}^{(Cp\log n)^{1/p}}%
\frac{\alpha(t)dt}{t^{1-p}}.
\]

For an important particular case (\ref{weib_like}), expressions for $a_{n}$
and $b_{n}$ will be evaluated below from (\ref{b_n}) by the same iteration
method, with using the second relation in (\ref{types}).

Using representation (\ref{expression2_gamma}), we get for all sufficiently
large $n$ that
\begin{align}
&  \gamma_{n}(x)-x=-\log\frac{1-F(b_{n}+a_{n}x)}{1-F(b_{n})}-x\nonumber\\
&  =\int_{b_{n}}^{b_{n}+f(b_{n})x}\frac{g(t)dt}{f(t)}-x=\int_{b_{n}}%
^{b_{n}+f(b_{n})x}\frac{1}{f(t)}dt+\int_{b_{n}}^{b_{n}+f(b_{n})x}\frac
{\alpha(t)dt}{f(t)}-x. \label{rate_weibull}%
\end{align}
Using the expression for $f$ and substituting $t=b_{n}v$, we get, that the
first integral on the right is equal to
\[
C^{-1}b_{n}^{p}\int_{1}^{1+Cb_{n}^{-p}x}\left(  v^{p-1}-1\right)  dv+x.
\]
Integrating and using Taylor for the integral (taking in mind that
$b_{n}\rightarrow\infty$ and $p>0$), we get that this is equal to
\[
\frac{1}{2}C(p-1)b_{n}^{-p}x^{2}(1+O(b_{n}^{-p}))=\frac{(p-1)x^{2}%
(1+o(1))}{2p\log n},\ \text{ }n\rightarrow\infty.
\]
Notice that for $p=1$ this integral vanishes. For the second integral on the
right in (\ref{rate_weibull}), changing $t=b_{n}+f(b_{n})v,$ we get that
\begin{equation}
\int_{b_{n}}^{b_{n}+f(b_{n})x}\frac{\alpha(t)dt}{f(t)}=(1+O(b_{n}^{-p}%
))\int_{0}^{x}\alpha(b_{n}+Cb_{n}^{1-p}v)dv,\ \ \text{ }n\rightarrow\infty.
\label{alpha}%
\end{equation}
Therefore
\[
\gamma_{n}(x)-x=\frac{(1+O(b_{n}^{-p}))}{2}C(p-1)b_{n}^{-p}x^{2}%
+(1+O(b_{n}^{-p}))\int_{0}^{x}\alpha(b_{n}+Cb_{n}^{1-p}v)dv
\]%
\begin{equation}
=(1+o(1))\left(  \frac{(p-1)x^{2}}{2p\log n}+\int_{0}^{x}\alpha(b_{n}%
+Cb_{n}^{1-p}v)dv\right)  ,\ \ n\rightarrow\infty. \label{rate_weibull1}%
\end{equation}
Thus we have proven the following refinement of Corollary
\ref{rate of convergence} for generalized Webull-like distributions.

\begin{proposition}
\label{rate_weibull-like} Let $F$ satisfies (\ref{weib_like}). Then for the
correction term in Corollary \ref{rate of convergence}, relation
(\ref{rate_weibull1}) is valid.
\end{proposition}

Obviously only for exponential like tail, $p=1$, the rate of convergence may
be better than logarithmic, it depends on the behavior of $\alpha(t)=g(t)-1.$
As we have already seen, if $\alpha(t)\equiv0,$ then, by virtue of Theorem
\ref{expansion}, the rate of convergence is proportional to $n^{-1}.$

\subsubsection{Weibull like distributions}

Now turn to classical Weibull like distributions. A distribution on
$\mathbb{R}_{+}$ with distribution function $F$ such that
\begin{equation}
1-F(x)=\mathbf{I}_{\{x\geq0\}}\ell(x)x^{\alpha}e^{-cx^{p}},\ \ \text{with
}p,\ c>0,\text{ }\alpha\in\mathbb{R}, \label{classicalWl}%
\end{equation}
and $\ell(x),$ slowly varying at infinity function, is called a Weibull-like
distribution. By Theorem 1.3.1, \cite{BGT}, for any slowly varying $\ell(x),$
some $x_{0}$ and all $x\geq x_{0}$,
\[
\ell(x)=c(x)\exp\left(  \int_{x_{0}}^{x}\delta(t)/tdt\right)  ,
\]
where $c(x)$ tends to a positive limit as $x\rightarrow\infty$ and
$\delta(t)\rightarrow0$ as $t\rightarrow\infty.$ Assume for simplicity that
$\ell(x)$ is normalized, \cite{BGT}, this means that $c(x)$ is a constant,
$c(x)=c>0.$ Hence, for some other positive $c$,
\begin{equation}
1-F(x)=c\exp\left(  \int_{x_{0}}^{x}\frac{\alpha+\delta(t)-pct^{p}}%
{t}dt\right)  . \label{alpha(t)0}%
\end{equation}
Hence in the model (\ref{weib_like}, \ref{weib_like1}),
\begin{equation}
C=\frac{1}{cp}\text{ \ and \ }\alpha(t)=-\frac{\alpha+\delta(t)}{cpt^{p}}.
\label{alpha(t)}%
\end{equation}
Inputting this in (\ref{rate_weibull1}) and integrating, we get after pretty
tedious but standard evaluations
\begin{equation}
\gamma_{n}(x)-x=\left(  \frac{(p-1)x^{2}}{2}-\alpha x+o(1)\right)  \frac
{1}{p\log n}\text{ } \label{Weibull-like rate 1}%
\end{equation}
as $n\rightarrow\infty.$

In order to evaluate normalizing sequences $b_{n}$ and $a_{n}$ on can apply
again the iteration method which has been applied in proof of Proposition
\ref{rate_weibull-like}. After pretty complicated but standard calculations,
see Appendix, we get the following.

\emph{For }$p=1$ \emph{one can take}
\begin{equation}
a_{n}=\frac{1}{c},\ b_{n}=\frac{1}{c}\log n+\frac{\alpha}{c}\log\left(
\frac{1}{c}\log n\right)  . \label{norm_p=1}%
\end{equation}

\emph{For } $p\neq1$ \emph{one can take}
\begin{equation}
a_{n}=\frac{1}{cp}\left(  \frac{1}{c}\log n\right)  ^{1/p-1};
\label{a_n_weibl}%
\end{equation}
\begin{align}
b_{n}  &  =\left(  c^{-1}\log n\right)  ^{1/p}\nonumber\\
&  +\frac{1}{p}\left(  c^{-1}\log n\right)  ^{1/p-1}\left(  \frac{a}{pc}%
\log(c^{-1}\log n)-c^{-1}\log\ell\left(  (c^{-1}\log n)^{1/p}\right)  \right)
. \label{b_n_weibl}%
\end{align}

Notice that in \cite{Gasull}, the sequences $a_{n}$ and $b_{n}$ has been
evaluates by another way, using so called $\mathbf{W}$ Lamperti functions, but
only in case $\ell(x)$ is constant. Thus we have the following refinement of
Corollary \ref{rate of convergence} for Webull-like distributions.

\begin{proposition}
\label{weibull_like_standad_rate} Let a distribution function $F$ be
ultimately absolutely continuous and (\ref{classicalWl}) takes place for it.
Then the correction term of Corollary \ref{rate of convergence} satisfies
(\ref{Weibull-like rate 1}). Moreover, one can take $a_{n}$ and $b_{n}$ as in
(\ref{a_n_weibl}), (\ref{b_n_weibl}), correspondingly.
\end{proposition}

\subsubsection{Example. Weibull distribution.}

Consider Wibull distribution, since it is extremely important in many fields,
such as reliability theory, queuing theory, finances. That is, let $\alpha=0,$
$\ell(x)\equiv1$ in (\ref{classicalWl}). One can take
\[
b_{n}=\left(  \frac{1}{c}\log n\right)  ^{1/p}+\frac{1}{pc}\log\frac{1}%
{c}\left(  \frac{1}{c}\log n\right)  ^{1/p-1},
\]
leaving the same $a_{n}$. After integrating, we get that
\begin{align*}
\gamma_{n}(x)-x  &  =\log n\left(  \left(  1+\frac{x}{p\log n}\right)
^{p}-1\right)  -x\\
&  =\frac{1}{\log n}\sum_{k=0}^{\infty}\binom{p}{k+2}\frac{x^{k+2}}%
{p^{k+2}\log^{k}n},\ \ \text{for }p\neq1;
\end{align*}
and
\[
\gamma_{n}(x)=x,\ \ \text{for }p=1.\
\]
Hence, if $p\neq1$, the decomposition in powers of $\log n$ takes place,
\[
G_{n}(x)=\exp\left(  -e^{-\gamma_{n}(x)}\right)  +O(n^{-1}),\ \ \text{\ \ }%
n\rightarrow\infty;
\]
for $p=1$, we have the decomposition in powers of $n$,
\[
G_{n}(x)=\exp\left(  -e^{-x}\left(  1+\sum_{k=1}^{\infty}\frac{(-1)^{k}%
e^{-kx}}{(k+1)n^{k}}\right)  \right)  .
\]
Remark that for $p=2,$
\[
\gamma_{n}(x)-x=\frac{x^{2}}{4\log n},
\]
which corresponds to the member $x^{2}/2l_{T}^{2}$ in (\ref{A_T}), with change
$n$ on $T/(2\pi).$ Remind that for Gaussian process $X(t)$ in Theorem
\ref{main}, for any $a$ and some $\delta>0$ we have,
\[
P\left(  \max_{t\in\lbrack0,a]}X(t)>u\right)  =\frac{a}{2\pi}e^{-u^{2}%
/2}+P(X(0)>u)+O(e^{-(1+\delta)u^{2}/2})
\]
as $u\rightarrow\infty,$ see \cite{lect}, where for proof of the Theorem we
take $a=a(T)\rightarrow\infty,$ but $a(T)/T\rightarrow0,$ $T\rightarrow\infty$
in order to make negligible the probability in the right hand part and to use
Weibull distribution with $p=2.$

\subsection{Generalized log-Weibull-like distributions}

Assume again that $1-F$ is ultimately absolutely continuous and take in
(\ref{fon MisesG}), $f(t)=Ct\log^{1-p}t$, $p>1.$ Also assume that conditions
(\ref{Fsmooth}) are fulfilled. Then redefining $g(t)$ on a finite interval, as
earlier, consider distribution functions satisfying
\begin{equation}
1-F(x)=\exp\left(  -\int_{1}^{x}\frac{g(t)dt}{Ct\log^{1-p}t}\right)
\label{log-weibG}%
\end{equation}
$x\geq1,$ with an agreement that $1/\log^{1-p}1=0.$ Remind that $\alpha
(t)=g(t)-1,$ and $\alpha(t)\rightarrow0$ asè $t\rightarrow\infty.$
Notice that for $p\leq1$ such the functions do not belong to $MDA(\Lambda)$.

Proof of the following refinement of Corollary \ref{rate of convergence} for
\emph{generalized log-Weibull-like distributions }(\ref{log-weibG}) is given
in Appendix.

\begin{proposition}
\label{log-Weib_prop} For generalized log-Weibull-like distributions
(\ref{log-weibG}), assertion of Corollary  \ref{rate of convergence} takes place with
\begin{equation}
\gamma_{n}(x)-x=(1+o(1))\left(  2^{-1}C^{1/p}p^{(1-p)/p}x^{2}\log
^{1/p-1}n+\int_{0}^{x}\alpha(b_{n}+Cb_{n}^{1-p}v)dv\right)
\label{log-Weib_prop-f}%
\end{equation}
as $n\rightarrow\infty$.
\end{proposition}

That is, by Corollary \ref{rate of convergence}, the rate of convergence in
(\ref{Gumbel_limit}) for such distributions is proportional to the right hand
part of this equality.

Remark that in the case of log-Weibull-like distributions, since $p>1$, the
rate of convergence cannot be better than logarithmic and also depends through
the second term above of how fast $\alpha(t)$ tends to zero as $t\rightarrow
\infty$. Remark also that for $p\leq1$ the distributions (\ref{log-weibG}) do
not belong to $MDA(\Lambda).$

\subsubsection{ Log-Weibull-like distributions}

Now consider a particular case of distributions (\ref{log-weibG}), A
distribution on $[1,\infty)$ with distribution function $F(x)$ such that
\begin{equation}
1-F(x)=\ell(x)x^{\alpha}e^{-c\log^{p}x},\ \ \text{with }p>1,\ c>0,\text{
}\alpha\in\mathbb{R}, \label{log-weib}%
\end{equation}
and $\ell(x),$ slowly varying on infinity function, is called a
\emph{log-Weibull-like distribution}. Similarly to representation
(\ref{alpha(t)0}), using (\ref{alpha(t)}) with corresponding modification
$\delta(t)$ on a finite interval, we get by simple calculus, that
\[
1-F(x)=\exp\left(  -\int_{1}^{x}\frac{cp}{t\log^{1-p}t}\left(  1-\frac
{\alpha+\delta(t)}{cp\log^{p-1}t}\right)  \right)  dt,
\]
so that
\begin{equation}
C=\frac{1}{cp},\text{ and }g(t)=1-\frac{\alpha+\delta(t)}{cp\log^{p-1}t}.
\label{log-Weib-alpha}%
\end{equation}
Now apply Proposition \ref{log-Weib_prop}. After simple calculations we get
from (\ref{b_nl}), (\ref{log-Weib-alpha}), that for $\alpha\neq0,$
\[
\alpha(b_{n}v)=-\frac{\alpha(1+o(1))}{cp\log^{p-1}(b_{n}v)}=-\frac
{\alpha(1+o(1))}{(cp)^{(2p-1)/p}\ln^{(p-1)/p}n},
\]
that is, the second member on the right in (\ref{log-Weib_prop-f}) represents
the rate of convergence. For $\alpha=0$ the rate can be better, in dependence
of rate of tending $\delta(t)$ to zero. In case $\ell(t)$ is constant, the
first member on the right in brakets of (\ref{log-Weib_prop-f}) is the rate of convergence.

Using again asymptotical iterations, one can evaluate expressions for $a_{n}$
and $b_{n},$ it is similar to corresponding evaluations for Weibull-like tails.

\section{ A scale in $MDA(\Lambda)$ for distributions with smooth tails}

Classes of generalized Weibull-like and log-Weibull-like distributions can be
a beginning of a natural scale in $MDA(\Lambda).$ As it is mentioned in
Introduction, this domain is enormous wide, distributions with very different
tail behaviors from this domain plays important role in financial, actuarial
research, reliability theory, engeneering aplications. Therefore Gnedenko
limit theorem (\ref{Gumbel_limit}) gives too far from complete information on
the tail behavior of a ditribution under estimation. For example, assignment
of insurance premium is strongly depends on behavior of tail distribution of
the time to insurance case. There are a plenty of studies in statistical
discrimination of hypotheses about the distribution tails for Weibull and
log-Weibull distributions, see \cite{Clerjaud}, \cite{rodionov},
\cite{rodionov1} and references therein.

A continuation of a scale which begins with the two considered here
distribution classes can be as following. In the two classes we have $f(t)=$
$Ct^{1-p},$ $p>0$ and $f(t)=Ct\log^{1-p}t,$ $p>1,$correspondingly.
Disgtibutions with heavier tails, say, with tails $\exp(-C\log x\log\log
^{p}x),$ $p>1$ can be described by von Mises function with $f(t)=Ct(\log\log
t)^{a},$ $a<0,$the constant $C$ can be different. The next "scale division" is
$f(t)=Ct(\log\log\log t)^{-a},$ so on. Obviously, for any natural $k,$
$f(t)=Ct(\log_{(k)}t)^{-a}$ satisfies conditions for representations
(\ref{fon Mises0}) and \ref{fon MisesG}). The repetitions number $k$ of
logarithms can be called Gumbel index, then the considered here classes of
distributions have indexex $k=0$ and $k=1,$ correspondingly.

It is interesting to consider behavior of $\gamma_{n}(x)-x$ as $n\rightarrow
\infty$ for all the scale. As we see, it is sufficient to consider only
representation (\ref{fon Mises0}), that is, $g(t)=1.$

\section{Appendix.}

\subsection{Derivation of normalizations (\ref{norm_p=1}, \ref{a_n_weibl},
\ref{b_n_weibl})}

Taking logarithm of
\[
\ell(x)x^{\alpha}e^{-cx^{p}}=n^{-1}%
\]
we have the equation for $b_{n},$
\[
\frac{1}{c}\log\ell(x)+\frac{\alpha}{c}\log x-x^{p}=-\frac{1}{c}\log n.
\]
Substituting $u=u_{n}=c^{-1}\log n,$ $y=x^{p},$ we have,
\begin{equation}
y-\frac{\alpha}{pc}\log y-\frac{1}{c}\log\ell_{1}(y)=u, \label{eq1}%
\end{equation}
where $\ell_{1}(y)=\ell(y^{1/p})$ is also slowly varying with $\delta
_{1}(t)=p\delta(t^{p})$ in given above representation for smooth slowly
varying functions. Now find asymptotically a root of (\ref{eq1}) as
$u\rightarrow\infty$ by asymptotic iteration. First, $y=u(1+o(1)).$
Substituting again this into the equation, we easily get, that
\[
y=u+L(u)+o(1),\ u\rightarrow\infty,
\]
where
\[
L(u):=\frac{\alpha}{pc}\log u+\frac{1}{c}\log\ell_{2}(u),
\]
and $\ell_{2}(y)=\ell_{1}(y(1+o(1)),$ $y\rightarrow\infty,$ is slowly varying
with corresponding $\delta_{2}(t),$ obvious modification of $\delta_{1}(t).$
Again inputting this in (\ref{eq1}), we get, that
\begin{align*}
y  &  =u+\frac{\alpha}{pc}\log\left(  u+L(u)+o(1)\right)  +\frac{1}{c}\log
\ell_{1}\left(  u+L(u)+o(1)\right) \\
&  =u+\frac{\alpha}{pc}\log u+\frac{\alpha}{pc}\log\left(  1+\frac
{L(u)+o(1)}{u}\right)  +\frac{1}{c}\log\ell_{1}\left(  u+L(u)+o(1)\right)
\end{align*}%
\[
=u+\frac{\alpha}{pc}\log u+\frac{\alpha}{pc}\frac{L(u)}{u}(1+o(1))+\frac{1}%
{c}\log\ell_{1}\left(  u+L(u)+o(1)\right)  ,\ u\rightarrow\infty.
\]
Consider the last summand on the right. It is equal to
\[
\int_{1}^{u+L(u)+o(1)}\frac{\delta_{1}(t)}{t}dt=\log\ell_{1}(u)+\Delta(u),
\]
where
\[
\Delta(u)=\int_{u}^{u+L(u)+o(1)}\frac{\delta_{1}(t)}{t}dt=\frac{\alpha}%
{pc}\frac{L(u)}{u}\frac{\delta_{1}(\theta u)}{\theta}%
\]
with%
\[
\theta\in\left[  1,1+\frac{L(u)+o(1)}{u}\right]  .
\]
Thus we have the following asymptotical equality,
\begin{align}
b_{n}^{p}  &  =u_{n}+L(u_{n})+\frac{\alpha}{pc}\frac{L(u_{n})}{u_{n}}%
+\frac{o(1)}{u_{n}}+\frac{1}{c}\Delta(u_{n})\nonumber\\
&  =u_{n}+\frac{\alpha}{pc}\log u_{n}+\frac{1}{c}\log\ell_{2}(u_{n})+\frac
{1}{u_{n}}R(u_{n}),\ \ n\rightarrow\infty. \label{R}%
\end{align}
It is easily follows from the above evaluations, that $R(u_{n})=O(L(u_{n}))$
as $n\rightarrow\infty.$

Furthermore,
\[
a_{n}=f(b_{n})=\frac{1}{cp}b_{n}^{1-p}.\
\]
That is, for $p=1$ one can take normalized sequences (\ref{norm_p=1}).
Further, notice that from (\ref{R}) it foolows that any normalizing sequence
$\tilde{b}_{n}$ is equal to $b_{n}+o(1).$ Hence for $p\neq1$, applying in
(\ref{R}), we get that
\begin{align}
b_{n}  &  =u_{n}^{1/p}+\frac{1}{p}u_{n}^{1/p-1}\left(  \frac{\alpha}{pc}\log
u_{n}+\frac{1}{c}\log\ell_{2}(u_{n})+\frac{1}{u_{n}}R(u_{n})\right)
\label{b_n1}\\
&  \frac{1}{2p}\left(  \frac{1}{p}-1\right)  u_{n}^{1/p-2}\left(  \frac
{\alpha}{pc}\log u_{n}+\frac{1}{c}\log\ell_{2}(u_{n})+\frac{1}{u_{n}}%
R(u_{n})\right)  ^{2}+...,\nonumber
\end{align}
If $p>1,$ then
\[
a_{n}=\frac{1}{cp}u_{n}^{1/p-1}(1+o(1)),\ n\rightarrow\infty,
\]
and $a_{n}$ can be taken as in (\ref{a_n_weibl}). Using that $\ell_{2}%
(y)=\ell(y^{1/p})+O(\Delta(y)),$we see that $b_{n}$ can be taken as in
(\ref{b_n_weibl})$.$ The same is valid also for $p<1.$ Indeed, in gthis case
$a_{n}\rightarrow\infty$ as $n\rightarrow\infty,$ and even if $p<1/2,$ the
other members in (\ref{b_n1}) can tend to infinity with $n$, equalities
(\ref{types}) let us to remind the same members in (\ref{b_n1}) as in case
$p>1.$

\subsection{Proof of Proposition \ref{log-Weib_prop}.}

Similarly to evaluations rewrite equation (\ref{log-weibG}) for $b_{n}$ as
\[
\log^{p}x-\int_{1}^{x}\frac{p\alpha(t)dt}{t\log^{1-p}t}=Cp\log n.
\]
Denote $y=\log^{p}x,$ $u=Cp\log n,$ change $s=\log t,$ the above equation
takes the form
\begin{equation}
\label{iter1l}y-p\int_{0}^{y^{1/p}}\alpha(e^{s})s^{p-1}ds=u.
\end{equation}
From here we get for the solution $y=u(1+o(1)),$ as $n\rightarrow\infty.$
Write for convenience $y=u(1+\varepsilon(y,u))^{p}.$ Putting this to
(\ref{iter1l}), we have,
\begin{align}
y  &  =u+p\int_{0}^{u^{1/p}(1+\varepsilon(y,u)}\alpha(e^{s})s^{p-1}%
ds\nonumber\\
&  =u+p\int_{0}^{u^{1/p}}\alpha(e^{s})s^{p-1}ds+p\int_{u^{1/p}}^{u^{1/p}%
(1+\varepsilon(y,u))}\alpha(e^{s})s^{p-1}ds\nonumber\\
&  =u+p\int_{0}^{u^{1/p}}\alpha(e^{s})s^{p-1}ds+pu^{-1}\int_{1}^{1+\varepsilon
(y,u)}\alpha(\exp(u^{-1/p}v))v^{p-1}dv, \label{iter2l}%
\end{align}
where we changed $v=u^{1/p}s$ in the last integral. Since $\alpha(1)$ can be
chosen to be a positive constant, the last integral is of order $\varepsilon
(y,u).$ Now, applying asymptotical iteration method, in order to find
$\varepsilon(y,u),$ we input this again into the upper integration limit in
(\ref{iter1}), so on. Finally we get from (\ref{iter2l}) that%
\begin{equation}
\log^{p}b_{n}=Cp\log n+p\int_{0}^{(Cp\log n)^{1/p}}\alpha(e^{s})s^{p-1}%
ds+\frac{\varepsilon^{\prime}(n,b_{n})}{\log n}, \label{b_nl}%
\end{equation}
where%
\[
\varepsilon^{\prime}(n,b_{n})=\frac{1}{C}\int_{1}^{1+\varepsilon(\log^{p}%
b_{n},Cp)\log n)}\alpha(n^{C}v))v^{p-1}dv.
\]

Now evaluate a rate of convergence in the limit theorem. Integrating in
(\ref{rate_weibull}) with a new èíòåã $f$, we
get, that
\begin{align*}
\int_{b_{n}}^{b_{n}+f(b_{n})x}\frac{dt}{f(t)}dt  &  =\frac{1}{C}\int_{b_{n}%
}^{b_{n}+f(b_{n})x}\log^{p-1}td\log t=\left.  \frac{1}{pC}\log^{p}t\right\vert
_{b_{n}}^{b_{n}+f(b_{n})x}\\
&  =\frac{1}{pC}\left(  \log^{p}(b_{n}(1+Cx\log^{1-p}b_{n}))-\log^{p}%
b_{n}\right)
\end{align*}
Applying several times Taylor expansion, after tedious but obvious evaluations we
get, that the difference in latter brakets is equal to
\begin{align*}
&  Cpx+\frac{1}{2}C^{2}px^{2}\log^{1-p}b_{n}+\frac{p(p-1)}{2}C^{2}x^{2}%
\log^{-p}b_{n}+O(\log^{2-2p}b_{n})\\
&  =Cpx+\frac{1}{2}C^{2}px^{2}\log^{1-p}b_{n}+O(\log^{-p}b_{n}+\log
^{2-2p}b_{n}).
\end{align*}
Hence
\[
\int_{b_{n}}^{b_{n}+f(b_{n})x}\frac{dt}{f(t)}dt=x+\frac{1}{2}Cx^{2}\log
^{1-p}b_{n}(1+O(\log^{-1\wedge(p-1)}b_{n})).
\]
Furthermore, similarly to (\ref{alpha}), changing $t=b_{n}+f(b_{n})v,$ we get
the following,%
\[
\int_{b_{n}}^{b_{n}+f(b_{n})x}\frac{\alpha(t)dt}{f(t)}=(1+O(\log^{1-p}%
b_{n}))\int_{0}^{x}\alpha(b_{n}(1+C\log^{1-p}v)dv,\ \ \text{ }n\rightarrow
\infty.
\]
Therefore%

\begin{align*}
\gamma_{n}(x)-x  &  =\frac{1+O(\log^{-1\wedge(p-1)}b_{n})}{2}Cx^{2}\log
^{1-p}b_{n}\\
&  +(1+O(\log^{1-p}b_{n}))\int_{0}^{x}\alpha(b_{n}+Cb_{n}^{1-p}v)dv,
\end{align*}
which follows (\ref{log-Weib_prop-f}).

\end{document}